\documentclass[a4paper,12pt]{article}

\usepackage[latin1]{inputenc}
\usepackage[english]{babel}
\usepackage{amssymb}
\usepackage{amsmath}
\usepackage{latexsym}
\usepackage{amsthm}
\usepackage[pdftex]{graphicx}
\usepackage{pgf, tikz} 
\usepackage{hyperref,url}
\usepackage{pgfplots}  
\usepackage{mathtools}
\pgfplotsset{width=6.6cm,compat=1.7}
\usepackage{array}
\usepackage{fancyhdr}
\usepackage{color}
\usetikzlibrary{arrows}

\usetikzlibrary{cd}

\usetikzlibrary{decorations.pathreplacing}

\DeclareMathOperator*{\Des}{Des}

\DeclareMathOperator*{\des}{des}

\newcommand{\seqnum}[1]{\href{http://oeis.org/#1}{\underline{#1}}}

 \theoremstyle{plain}
 \newtheorem{thm}{Theorem}[section]
 \newtheorem{cor}[thm]{Corollary}
 \newtheorem{lem}[thm]{Lemma}
 \newtheorem{prop}[thm]{Proposition}
 \newtheorem{conj}[thm]{Conjecture}

 \theoremstyle{definition}
 
 \newtheorem{example}[thm]{Example}

 \theoremstyle{remark}
 
\title{Pattern avoiding alternating involutions}
\date{}
\author{Marilena Barnabei\\
Dipartimento di Matematica \\
Universit\`a di Bologna, 40126, ITALY \\
\texttt{marilena.barnabei@unibo.it}\and
Flavio Bonetti \\
P.A.M. \\
Universit\`a di Bologna, 40126, ITALY \\
\texttt{flavio.bonetti@unibo.it}\and
Niccol\`o Castronuovo \\
Liceo ``A. Einstein'', Rimini, 47923, ITALY \\
\texttt{niccolo.castronuovo2@unibo.it}\and
Matteo Silimbani \thanks{corresponding author} \\
Istituto Comprensivo ``E. Rosetti'', Forlimpopoli, 47034, ITALY \\
\texttt{matteosilimbani@icrosetti.istruzioneer.it}}

\begin{document}
\maketitle
\begin{abstract}
We enumerate and characterize some classes of alternating and reverse alternating involutions avoiding a single pattern of length three or four. If on one hand the case of patterns of length three is trivial, on the other hand the length four case is more challenging and involves sequences of combinatorial interest, such as Motzkin and Fibonacci numbers. 
\end{abstract}

2020 msc: 05A05, 05A15 (primary), 05A19 (secondary).


\section{Introduction}

A permutation $\pi$ avoids a pattern $\tau$ whenever $\pi$ does not contain any subsequence order-isomorphic to $\tau.$

The theory of permutation patterns goes back to the work of Knuth \cite{Kn2}, who, in the
1970's, introduced the definition of pattern avoidance in connection to the stack sorting problem. The first systematic study of these objects appears in the paper by Simion and Schmidt (\cite{Si}). Nowadays, the theory is very rich and widely expanded, with hundreds of papers appeared in the last decades (see e.g. \cite{Ki} and references therein). 

More recently permutation patterns have been studied over particular subset of the symmetric group. 
In particular pattern avoidance has been studied over involutions (see e.g. \cite{Ba6,Bona,Jaggard,Ki}) and over alternating permutations (see e.g. \cite{Bonaf,ChenChenZhou,Gowravaram,Ki,Lewis,Lewis2,Lewis3,XuYan,Yan}), i.e.,
permutations $\pi=\pi_1\ldots \pi_n$ such that $\pi_i<\pi_{i+1}$ if and only if $i$ is odd.

The enumeration of alternating involutions is due to Stanley (see \cite{Stanleyalt1} and also his survey on alternating permutations \cite{Stanleyalt2}). However, to the best of our knowledge, pattern avoiding alternating involutions have not been studied so far.

In this paper we consider alternating involutions that avoid some patterns of length three or four. If on one hand the case of patterns of length three is trivial, on the other hand the length four case is more challenging and involves sequences of combinatorial interest, such as Motzkin and Fibonacci numbers.

\section{Preliminaries}

\subsection{Permutations}

Let $S_n$ be the symmetric group over the symbols $\{1,2,\ldots \,n\}.$ We will often write permutations in $S_n$ in one line notation as $\pi=\pi_1\ldots \pi_n.$
An \textit{involution} is a permutation $\pi$ such that $\pi=\pi^{-1}.$
The \textit{reverse} of a permutation $\pi=\pi_1\pi_2\ldots \pi_n$ is $\pi^r=\pi_n\pi_{n-1}\ldots \pi_1$ and the \textit{complement} of $\pi$ is $\pi^c=n+1-\pi_1\,n+1-\pi_2\ldots n+1-\pi_n.$ The \textit{reverse-complement} of $\pi$ is $\pi^{rc}=(\pi^r)^c=(\pi^c)^r.$
Notice that $(\pi^{-1})^{rc}=(\pi^{rc})^{-1}.$ In particular the reverse-complement of an involution is an involution.

A \textit{descent} in a permutation $\pi$ is an index $i$ such that $\pi_i>\pi_{i+1}.$
Denote by $\Des(\pi)$ the set of descents of $\pi$ and by $\des(\pi)$ its cardinality. 

Recall that a \textit{left-to-right minimum} (ltr minimum from now on) of a permutation $\pi$ is a value $\pi(i)$ such that $\pi(j)>\pi(i)$ for every $j<i.$ 
A \textit{left-to-right maximum} (ltr maximum from now on) is a value $\pi(i)$ such that $\pi(j)<\pi(i)$ for every $j<i.$ 
The definition of right-to-left (rtl) minimum and maximum is analogous.

\begin{lem}\label{mininv}
In any involution $\pi$ the ltr minima  form an involution themselves, i.e., a ltr minimum $m$ is at position $i$ if and only if $i$ is a ltr minimum. The same is true for rtl maxima.
\end{lem}
\proof
We will consider only the case of ltr minima, since the other case is analogous.
Denote by $m_1, m_2, \ldots, m_k$ the left-to-right minima of $\pi$, and by $i_1,i_2\ldots,i_k$ their respective positions. We want to show that $$\{m_1, m_2, \ldots, m_k\}=\{i_1,i_2\ldots,i_k\}.$$
Suppose on the contrary that there exists an index $j$ such that $i_j$ is not
the value of a left-to-right minimum. Then, to the left of $i_j$ in $\pi$ there is a symbol $a$ less than $i_j.$ In other terms, there exist two integers $a, h$ such that $\pi(h)=a,$ $a<i_j$ and $h<m_j.$ In this situation, $m_j$ is not a left-to-right minimum, since it is preceded by $h,$ and $h<m_j,$ contradiction.
\endproof

\begin{example}
Consider the involution $\pi=7\,9\,4\,3\,5\,6\,1\,10\,2\,8.$ The ltr minima of $\pi$ are 7,4,3,1 and they form an involution themselves. The same is true for the rtl maxima 8 and 10. 
\end{example}

Given a word $w=w_1\ldots w_j$ whose letters are distinct numbers, the \textit{standardization} of $w$ is the unique permutation $\pi$ in $S_j$ order isomorphic to $w.$ If two words $w$ and $u$ have the same standardization we write $u\sim w.$

The \textit{decomposition into connected components} of a permutation $\pi$ is the finest way to write $\pi$ as $\pi=w_1w_2\ldots w_k,$ where each $w_i$ is a permutation of the symbols from $|w_1|+|w_2|+\ldots+|w_{i-1}|+1$ to $|w_1|+|w_2|+\ldots+|w_{i-1}|+|w_i|.$ Each $w_i$ is called \textit{a connected component} of $\pi.$ A permutation is said to be \textit{connected} if it is composed by only one connected component. 

\begin{example}
The decomposition into connected components of the permutation $\pi=34215786$ is $\pi=w_1w_2w_3,$ where $w_1=3421,$ $w_2=5,$ $w_3=786.$
\end{example}

\subsection{Alternating permutations}

A permutation $\pi=\pi_1\ldots \pi_n$ is said to be \textit{alternating} if $\pi_i<\pi_{i+1}$ if and only if $i$ is odd and \textit{reverse alternating} if and only if $\pi_i>\pi_{i+1}$ if and only if $i$ is odd. Equivalently, a permutation $\pi$ is alternating whenever $\Des(\pi)=\{2,4,6,\ldots\}.$

\begin{example}
The permutation $\pi=4615273$ is alternating, while $\sigma=5372614$ is reverse alternating.
\end{example}

Denote by $S_n$ ($I_n,$ $A_n,$ $RA_n,$ $AI_n$ and $RAI_n,$ respectively) the set of permutations  (involutions, alternating permutations, reverse alternating permutations, alternating involutions, reverse alternating involutions, respectively) of length $n.$

The following lemma will be useful in the sequel.

\begin{lem}\label{str}
\begin{itemize}
\item[i)] Let $\pi=\pi_1\ldots\pi_n \in RAI_{n}.$ Then $1$ is in even position, $n$ is in odd position, hence
 $\pi_1$ is even and $\pi_{n}$ is odd.
\item[ii)] Let $\pi=\pi_1\ldots\pi_n \in AI_n.$ Then $1$ is in odd position, $n$ is in even position, hence
$\pi_1$ is odd and $\pi_n$ is even. 
\end{itemize}
\end{lem}

\subsection{Pattern avoidance}

A permutation $\sigma=\sigma_1\ldots\sigma_n \in S_n$ \textit{avoids} the  pattern $\tau \in S_k$ if there are no indices $i_1,i_2,\ldots,i_k$ such that the subsequence
$\sigma_{i_1}\sigma_{i_2}\ldots \sigma_{i_k}$ is order isomorphic to $\tau.$

Denote by $S_n(\tau)$ the set of permutations of length $n$ avoiding $\tau$ and let
$S(\tau):=\bigcup_n S_n(\tau).$ We will keep this notation also when $S_n$ is replaced by other subsets of permutations, such as $I_n,$ $A_n,$ etc.

Notice that an involution avoids $\tau$ if and only if it avoids $\tau^{-1}.$

The next trivial lemma will be useful in the sequel.

\begin{lem}\label{rc}
 The reverse-complement map is a bijection between $AI_{2n+1}(\tau)$ and $RAI_{2n+1}(\tau^{rc}),$ between $AI_{2n}(\tau)$ and $AI_{2n}(\tau^{rc})$ and between $RAI_{2n}(\tau)$ and $RAI_{2n}(\tau^{rc}).$ 
\end{lem}

\subsection{Motzkin paths}

A \textit{Motzkin path} of length $n$ is a lattice path starting at $(0,0),$ ending at $(n,0),$ consisting of up steps $U$
of the form $(1,1),$ down steps $D$ of the form $(1,-1),$ and horizontal steps $H$ of the form $(1,0),$ and lying
weakly above the $x$-axis.
As usual, a Motzkin path can be identified with a \textit{Motzkin word}, namely, a word $w = w_1w_2\ldots w_n$ of
length $n$ in the alphabet $\{U, D, H\}$ with the constraint that the number of occurrences of the letter $U$ is
equal to the number of occurrences of the letter $D$ and, for every $i,$ the number of occurrences of $U$ in
the subword $w_1w_2\ldots w_i$
is not smaller than the number of occurrences of $D.$ In the following we will not distinguish between a Motzkin path and the corresponding word. Denote by $\mathcal{M}_n$ the set of Motzkin path of length $n$ and by $M_n$ its cardinality, the $n$-th Motzkin number (see sequence \seqnum{A001006} in \cite{Sl}).
The \textit{diod decomposition} of a Motzkin path $m$ of even length $2n$ is the decomposition of $m$ as $m=d_1d_2\ldots d_n,$ where each $d_i$ is a subword of $m$ of length two. Each $d_i$ in the diod decomposition of a Motzkin path is called a \textit{diod} of $m.$

\section{General results}

In this section we prove two general results which will be used in the paper. 

\begin{lem}\label{BWX}
Let $\tau$ be any permutation of $\{3,\ldots ,m\},$ $m\geq 3,$ Then $$|AI_{n}(12\tau)|=|AI_n(21\tau)|.$$

\end{lem}
\proof
We closely follow  \cite{Ouc}, where a similar result is proved for doubly alternating permutations, i.e., alternating permutations whose inverse is also alternating. 

Our goal is to find a bijection between $AI_{n}(12\tau)=AI_n(12\tau,12\tau^{-1})$ and $AI_{n}(21\tau)=AI_n(21\tau,21\tau^{-1}).$ 

We will use the diagrammatic representation of a permutation, i.e., we will identify a permutation $\pi$ in $S_n$ with a $n\times n$ diagram with a dot at position $(i,\pi_i),$ for every $i$
(notice that Ouchterlony \cite{Ouc} uses a slightly different definition for the diagram of a permutation).

A dot, $d,$ in the diagram of a permutation $\pi$ is called \textit{active} if $d$ is the 1 or 2 in any $12\tau,$ $12\tau^{-1},$ $21\tau$ or $21\tau^{-1}$ pattern in $\pi,$ or \textit{inactive} otherwise.
Also the pair of dots $(d_1,d_2)$ is called an \textit{active pair} if $d_1 d_2$ is the $12$ in a $12\tau$ or $12\tau^{-1}$ pattern or the $21$
in a $21\tau$ or $21\tau^{-1}$ pattern.


We now  define a Young diagram, $\lambda_{\pi},$ consisting of the part of the diagram of $\pi$ which contains the active dots.
For any two dots $d_1,d_2,$ let $R_{d_1,d_2}$  be the smallest rectangle with bottom left coordinates (1, 1), such that $d_1,d_2\in R_{d_1,d_2}.$ Define $$\lambda_{\pi}=\bigcup R_{d_1,d_2},$$
where the union is over all active pairs $(d_1,d_2)$ of $\pi.$ It is clear from the definition that $\lambda_{\pi}$ is indeed a Young
diagram. 
Since $\pi$ is an involution, its diagram is symmetric with respect to the main diagonal, and for every active dot in position $(i,j)$ there is also an active dot in position $(j,i),$ hence $\lambda_{\pi}$ is a Young diagram symmetric with respect to the main diagonal. 

A \textit{rook placement} of a Young diagram $\lambda$ is a placement of dots in its boxes, such that all rows and columns contain exactly one dot. If some of the rows or columns are empty we call
it a \textit{partial rook placement}. Furthermore, we say that a rook placement on $\lambda$ avoids the pattern $\tau$ if no
rectangle, $R\subseteq \lambda,$ contains $\tau.$ 
Notice that the rook placement on $\lambda_\pi$ induced by an involution $\pi$ is symmetric.


Now, it has been proved by Jaggard \cite[Theorem 4.2]{Jaggard}, that the number of symmetric rook placements on the self-conjugate shape $\mu$ avoiding the patterns $12\tau$ and $12\tau^{-1}$ is equal to the number of 
symmetric rook placements on the self-conjugate shape $\mu$ avoiding the patterns $21\tau$ and $21\tau^{-1}.$

We call two permutations $\pi$ and $\sigma$ of the same size \textit{a-equivalent} if they have  the same inactive dots, and write $\pi\sim_a \sigma.$

In the sequel we will need the following three facts that are immediate consequences  of Lemma 6.4, Lemma 6.3 and Lemma 6.2 in \cite{Ouc}, respectively.
\begin{itemize}
\item[\textbf{(I)}] $\pi \sim_a \sigma \Rightarrow \lambda_\pi=\lambda_\sigma.$
\item[\textbf{(II)}] If $\pi\in AI_n(12\tau,12\tau^{-1})\cup AI_n(21\tau,21\tau^{-1})$ and $\pi\sim_a \sigma$ then $\sigma$ is doubly alternating.
\item[\textbf{(III)}] If $\pi\in AI_n$ and $rp(\lambda_\pi)$ is the partial rook placement on $\lambda_\pi$ induced by $\pi,$ then $$\pi \in AI_n(12\tau,12\tau^{-1}) \Leftrightarrow rp(\lambda_\pi) \mbox{ is }12\mbox{-avoiding}$$ and $$\pi \in AI_n(21\tau,21\tau^{-1}) \Leftrightarrow rp(\lambda_\pi) \mbox{ is }21\mbox{-avoiding}.$$
\end{itemize}


Now we are ready to construct a bijection $$\phi:AI_n(12\tau,12\tau^{-1})\to AI_n(21\tau,21\tau^{-1}).$$
Let $\pi \in AI_n(12\tau,12\tau^{-1}),$ so that the restriction of $\pi$ to $\lambda_\pi$ is a partial $12$-avoiding symmetric rook placement. By Jaggard's theorem (ignoring the empty rows and
columns) and by point \textbf{(I)} there exists a unique $21$-avoiding (partial) symmetric rook placement on $\lambda_{\pi}$ with the same rows and columns empty, which we combine with the inactive dots of $\pi$ to get $\phi(\pi).$ By point \textbf{(II)}, $\phi(\pi)$ is doubly alternating and, since its diagram is symmetric, it is an involution. By point \textbf{(III)} it avoids $21\tau$ (and hence also $21\tau^{-1}$). It is also clear from Jaggard's theorem that $\phi$ is indeed a bijection.
 
\endproof

\begin{example}
Consider the permutation $\pi=593716482\in AI_{9}(12435).$ Notice that $\pi$ contains the pattern $21435.$
Here $\tau=435.$
The diagram of $\pi$ is the following

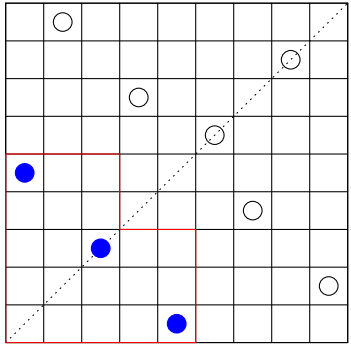
\centering
\begin{tikzpicture}[line cap=round,line join=round,>=triangle 45,x=0.5cm,y=0.5cm]
\draw[dotted,thin] (1,1)--(10,10);
\draw [] (1.,1.)-- (10.,1.);
\draw [] (10.,1.)-- (10.,10.);
\draw [] (10.,10.)-- (1.,10.);
\draw [] (1.,10.)-- (1.,1.);
\draw [] (9.,10.) -- (9.,1.);
\draw [] (10.,9.) -- (1.,9.);
\draw (2.,1.)-- (2.,10.);
\draw (3.,1.)-- (3.,10.);
\draw (4.,1.)-- (4.,10.);
\draw (5.,1.)-- (5.,10.);
\draw (6.,1.)-- (6.,10.);
\draw (7.,1.)-- (7.,10.);
\draw (8.,1.)-- (8.,10.);
\draw (1.,8.)-- (10.,8.);
\draw (1.,7.)-- (10.,7.);
\draw (1.,6.)-- (10.,6.);
\draw (1.,5.)-- (10.,5.);
\draw (1.,4.)-- (10.,4.);
\draw (1.,3.)-- (10.,3.);
\draw (1.,2.)-- (10.,2.);
\draw (1.,1.)-- (10.,1.);
\draw (1.,10.)-- (10.,10.);
\draw (1.,10.)-- (1.,1.);
\draw (10.,10.)-- (10.,1.);
\draw [color=red] (1.,6.) -- (4.,6.) -- (4.,4.) -- (6.,4.) -- (6.,1.);
\draw [color=red] (6.,1.) -- (1.,1.) -- (1.,6.);

\begin{scriptsize}
\draw [color=blue,fill=blue] (1.5,5.5) circle (3.5pt);
\draw [] (2.5,9.5) circle (3.5pt);
\draw [color=blue,fill=blue] (3.5,3.5) circle (3.5pt);
\draw [] (4.5,7.5) circle (3.5pt);
\draw [color=blue,fill=blue] (5.5,1.5) circle (3.5pt);
\draw [] (6.5,6.5) circle (3.5pt);
\draw [] (7.5,4.5) circle (3.5pt);
\draw [] (8.5,8.5) circle (3.5pt);
\draw [] (9.5,2.5) circle (3.5pt);

\end{scriptsize}
\end{tikzpicture}

\vspace{0.5cm}

The blue dots are the active dots of $\pi$ and the red Young diagram is $\lambda_{\pi}.$ Now applying to $\lambda_{\pi}$ the procedure described in the above proof we get $\phi(\pi)=195736482\in AI_9(21435).$
\end{example}

It follows from the previous lemma and Lemma \ref{rc}  that
if $\tau$ is any permutation of $\{1,\ldots ,k-2\}$ then $$|AI_{2n}(\tau\, k-1\, k)|=|AI_{2n}(\tau\, k\,k-1)|$$ and $$|RAI_{2n+1}(\tau\, k-1\, k)|=|RAI_{2n+1}(\tau\, k\,k-1)|.$$

Notice that similar relations do not hold for $RAI_{2n},$ in fact,  numerical computations show that, for general $n,$ $|RAI_{2n}(1234)|\neq |RAI_{2n}(1243)|$ and $|RAI_{2n}(1234)|\neq |RAI_{2n}(2134)|.$

When $\tau$ is an increasing sequence it is possible to provide a more explicit bijection $f$ between $|AI_{2n}(\tau\, k-1\, k)|$ and $|AI_{2n}(\tau\, k\,k-1)|.$ Such a bijection has been defined for the first time by J. West \cite{West} for permutations  and has been used by M. B\'ona \cite{Bonaf} to prove a conjecture by J. B. Lewis about alternating permutations. To this aim B\'ona proved that $f$ preserves the alternating property when the length of the permutation is even.
Here we recall the definition of the map $f:S_{t}(12\ldots \, k-1\, k)\to S_t(12\ldots\, k\,k-1).$  Consider a permutation $\pi \in S_{t}(12\ldots \, k-1\, k)$ and define the \textit{rank} of the element $\pi_i$ to be the maximum length of an increasing subsequence ending at $\pi_i.$ Since $\pi$ avoids $12\ldots \, k-1\, k,$ the maximal rank of an element is $k-1.$ Let $R$ be the set of  elements of $\pi$ whose rank is $k-1,$ and $P$  the set of their positions. 
The permutation $\rho=f(\pi)$ is obtained as follows.
\begin{itemize}
    \item if $j\notin P,$ $\rho_j=\pi_j,$
    \item if $j\in P,$ $\rho_j$ is the smallest unused element of $R$ that is larger than the closest entry of rank $k-2$ to the left of $\pi_j.$
\end{itemize}
Notice that, if $k=3,$ the map $f$ reduces to the classic Simion-Schmidt bijection (\cite{Si}, \cite{West}).

In the following lemma we prove that  $f$ preserves also the property of being an involution and we describe what happens in the odd case. 
\begin{lem}\label{Bonabij}
The map $f$ is a bijection between $AI_{2n}(12\ldots\, k-1\, k)$ and $AI_{2n}(12\ldots\, k\,k-1).$ Moreover, $f$ maps bijectively the subset of $AI_{2n+1}(1234)$ of those permutations having $2n+1$ at position 2 to the set $AI_{2n+1}(1243).$ 
\end{lem}
\proof
First of all we notice that, in an involution $\pi,$ all the elements of a given rank form an involution. To prove this fact it is sufficient to observe that elements of $\pi$ of rank 1 are the ltr minima of the permutation, which form an involution (see Lemma \ref{mininv}). We can consider the involution obtained removing from $\pi$ the elements of rank 1 and proceed inductively. 

Now we want to prove that $f$ sends involutions of lenght $n$ to involutions of length $n$ (for every length $n$). Also in this case we can obtain the result by induction on $k.$ When $k=3,$ we noticed above that $f$ coincides with the Simion-Schmidt map (\cite{Si}), and it is well-known that this map sends involutions to  involutions. If $k>3,$ we can delete from $\pi$ the elements of rank 1 and standardize, obtaining an involution $\pi'$ which avoids $12\ldots\, k-2\;\, k-1.$ If we apply the map $f,$ we obtain an involution avoiding $12\ldots\, k-1\;\, k-2$ which can be obtained from $f(\pi)$ by removing the elements of rank 1 and standardizing. Hence also $f(\pi)$ is an involution. The fact that $f$ preserves the alternating property (when the length is even) has been proved by B\'ona \cite[Theorem 1]{Bonaf}.

Now we turn our attention to the second assertion. We want to prove that $f^{-1}$ maps bijectively the set $AI_{2n+1}(1243)$ onto the set $Y_{2n+1}=AI_{2n+1}(1234)\cap \{\pi\in S_{2n+1}\,|\, \pi(2)=2n+1\}.$

The fact that the injective map $f^{-1}$ preserves the alternating property has been already observed by B\'ona \cite[Corollary 1]{Bonaf}. Moreover $f^{-1}$ preserves the involutory property, as proved above. 

Now consider $\pi\in AI_{2n+1}(1243).$ The symbol $2n+1$ appears at  even position $a.$ Suppose that $a\geq 4.$ Then $\pi_2=2n,$ otherwise $\pi_1\;\,\pi_2\;\,2n+1\;\,2n\sim 1243.$ As a consequence, $\pi_{2n+1}=a>2=\pi_{2n}$ which is impossible since $\pi$ is alternating. Hence $a=2$ and $\pi_2=2n+1.$

Since for every $\pi\in AI_{2n+1}(1243),$ $\pi(2)=2n+1,$  the symbols $2n+1$ and $2$ have rank $2.$ The map $f^{-1}$ fixes the elements of rank $k-2$ or less, hence it fixes the positions of $2n+1$ and $2.$

With similar arguments one can show that $f(Y_{2n+1})\subseteq AI_{2n+1}(1243)$ and conclude that $f(Y_{2n+1})= AI_{2n+1}(1243).$

\endproof

\begin{example}
Consider the permutation $\pi=5\,9\,7\,10\,1\,8\,3\,6\,2\,4\in AI_{10}(1234).$ The elements of rank 1 are 5 and 1, the elements of rank 2 are 9,7,3 and 2, the elements of rank 3 are 10, 8, 6 and 4. Then $f(\pi)=5\,9\,7\,8\,1\,10\,3\,4\,2\,6\in AI_{10}(1243).$
\end{example}

\section{Patterns of length three}\label{ltr}

The following proposition is an easy consequence of Proposition 3.1 in \cite{Ouc}. See also \cite{Oucth}.

\begin{prop}
$$|AI_n(123)|=|AI_n(213)|=|AI_n(231)|=|AI_n(231)|=$$
$$|RAI_n(132)|=|RAI_n(231)|=|RAI_n(312)|=|RAI_n(321)|=1$$
$$|AI_n(132)|=|RAI_n(213)|=\begin{cases} 1& \mbox{ if }n\mbox{ is even or }n=1\\
                            0& \mbox{ otherwise}.\end{cases} $$
$$|AI_n(321)|=|RAI_n(123)|=\begin{cases} 2& \mbox{ if }n\mbox{ is even and }n\geq 4\\
                            1& \mbox{ otherwise}.\end{cases} $$  
\end{prop}

\section{Patterns 4321 and 1234}\label{sec4321}

We recall the definition of a bijection $\Phi$ between the set $I_n(4321)$ of $4321$-avoiding involutions of the symmetric group $S_n$ and the set $\mathcal M_n$ of Motzkin paths of length $n$ (see \cite{Ba6}), that is essentially a restriction to the set $I_n(4321)$ of the bijection appearing in \cite{BIANE}.

Consider an involution $\tau$ avoiding $4321$ and determine the set $exc(\tau)=\{i|\tau(i) > i\}$ of its excedances. Start from the empty path and create the  Motzkin path $\Phi(\tau)$ by adding a step for every integer $1 \leq i \leq n$
as follows:
\begin{itemize}
    \item if $\tau(i)=i,$ add a horizontal step at the $i$-th position;
    \item if $\tau(i)>i,$ add an up step at the $i$-th position;
    \item if $\tau(i)<i,$ add a down step at the $i$-th position.
\end{itemize}

The map $\Phi$ is a bijection  whose inverse can be described as follows.

Given a Motzkin path $M,$ let $A=(a_1,\ldots,a_k)$ be the list  of  positions of up steps in $M,$ written in increasing order,  and let $B=(b_1,\ldots, b_k)$ the analogous list of  positions of  down steps. Then $\Phi^{-1}(M)$ is the involution $\tau$ given by the product of the cycles $(a_i,b_i)$ for $1\leq i\leq k.$

\begin{lem}
 $\Phi$ maps the set $RAI_{2n}(4321)$ onto the set $\widehat{\mathcal M}_{2n}$ of Motzkin paths of length $2n$ whose diods are either of the form $UH,$ $HD$ or $UD.$
\end{lem}
\proof
Let $\tau \in RAI_{2n}(4321).$ For every integer $k=1,2,\ldots,2n-1,$ if $k$ is odd, and hence a descent of $\tau$ , $\tau_k\tau_{k+1}$ is mapped by $\Phi$ to one of the following pairs: $UH,$ $UD,$ $HD.$
\endproof

\begin{thm}\label{RAI4321even}
$$|RAI_{2n}(4321)|=M_n.$$
\end{thm}
\proof
The previous lemma allows us to define a bijection $\hat{\Phi}:RAI_{2n}(4321)\to \mathcal{M}_n.$
Let $m=d_1d_2\ldots d_{n}$ be a Motzkin path (of even length) whose diods $d_i$ are either of the form $UH,$ $UD$ or $HD.$ 
Define $\Delta(m):=s_1\ldots s_n\in \mathcal M_n$  where
$$s_i=\begin{cases} U & \mbox{ if }d_i=UH\\   D & \mbox{ if }d_i=HD\\ H & \mbox{ if }d_i=UD.
\end{cases}$$
Set now $\hat{\Phi}=\Delta \circ \Phi.$
\endproof

\begin{example}
Consider $\pi=6\,2\,8\,4\,5\,1\,10\,3\,9\,7\in RAI_{10}(4321).$

Then $$\Phi(\pi)=\begin{tikzpicture}[every node/.style={draw,shape=circle,minimum size=1mm, inner sep=0mm, fill=black}]
\node[draw,shape=circle,minimum size=1mm, inner sep=0mm, fill=black] (10) at (10/2,0) {};
\node[draw,shape=circle,minimum size=1mm, inner sep=0mm, fill=black] (9) at (9/2,1/2) {};
\node[draw,shape=circle,minimum size=1mm, inner sep=0mm, fill=black] (8) at (8/2,1/2) {};
\node[draw,shape=circle,minimum size=1mm, inner sep=0mm, fill=black] (7) at (7/2,2/2) {};
\node[draw,shape=circle,minimum size=1mm, inner sep=0mm, fill=black] (6) at (6/2,1/2) {};
   \node[draw,shape=circle,minimum size=1mm, inner sep=0mm, fill=black] (5) at (5/2,1) {};
   \node[draw,shape=circle,minimum size=1mm, inner sep=0mm, fill=black] (4) at (4/2,2/2) {};
   \node[draw,shape=circle,minimum size=1mm, inner sep=0mm, fill=black] (3) at (3/2,2/2) {};
   \node[draw,shape=circle,minimum size=1mm, inner sep=0mm, fill=black] (2) at (2/2,1/2) {};
   \node[draw,shape=circle,minimum size=1mm, inner sep=0mm, fill=black] (1) at (1/2,1/2) {};
   \node[draw,shape=circle,minimum size=1mm, inner sep=0mm, fill=black] (0) at (0,0) {};
   \draw[] (0) -- (1) -- (2) -- (3) --   (4) -- (5) -- (6) -- (7) --(8)--(9)--(10);
\end{tikzpicture} $$

and $$\widehat{\Phi}(\pi)=\begin{tikzpicture}[every node/.style={draw,shape=circle,minimum size=1mm, inner sep=0mm, fill=black}]
   \node[draw,shape=circle,minimum size=1mm, inner sep=0mm, fill=black] (5) at (5/2,0) {};
   \node[draw,shape=circle,minimum size=1mm, inner sep=0mm, fill=black] (4) at (4/2,1/2) {};
   \node[draw,shape=circle,minimum size=1mm, inner sep=0mm, fill=black] (3) at (3/2,1/2) {};
   \node[draw,shape=circle,minimum size=1mm, inner sep=0mm, fill=black] (2) at (2/2,2/2) {};
   \node[draw,shape=circle,minimum size=1mm, inner sep=0mm, fill=black] (1) at (1/2,1/2) {};
   \node[draw,shape=circle,minimum size=1mm, inner sep=0mm, fill=black] (0) at (0,0) {};
   \draw[] (0) -- (1) -- (2) -- (3) --   (4) -- (5);
\end{tikzpicture}$$

\end{example}

\begin{thm}\label{RAI4321}
$$|RAI_{2n-1}(4321)|=M_n-M_{n-2}.$$
\end{thm}
\proof
Consider a permutation $\pi \in RAI_{2n}(4321).$ By Lemma \ref{str}, $\pi_{2n}$ is odd. Now consider the element $2n-1$ in $\pi.$ Notice that this element is at position $2n$ or $2n-1,$ otherwise $\pi$ would contain the pattern $2n\;\, 2n-1\;\, \pi(2n-1)\;\, \pi(2n) \sim 4321$ (notice that $2n$ appears before $2n-1$ because $\pi$ is an involution and $\pi_{2n-1}>\pi_{2n}$).

If $\pi_{2n}=2n-1,$ the permutation $\pi'$ obtained from 
$\pi$ by removing the element $2n-1$ and standardizing is an arbitrary element of $RAI_{2n-1}(4321)$ which ends with its maximum element.

If $\pi_{2n-1}=2n-1,$ write $\pi$ as 
$\pi=\left(\begin{smallmatrix} 
\cdots & k & \cdots & 2n-2 & 2n-1 & 2n\\
\cdots & 2n & \cdots & j & 2n-1 & k
\end{smallmatrix}\right)$.

If $k>j,$ consider the permutation $\pi'$ obtained from $\pi$  by removing the element $2n-1$ and standardizing. $\pi'$ is an arbitrary element of $RAI_{2n-1}(4321)$ which does not end with its maximum element. 

If $k<j,$ we must have $k=2n-3$ and $j=2n-2,$ otherwise $\pi$ would contain $4321.$ Hence we can write $\pi=\pi'\;\,2n\;\,2n-2\;\,2n-1\;\,2n-3,$ where $\pi'$ is an arbitrary permutation in $RAI_{2n-4}(4321).$

As a consequence, $|RAI_{2n}(4321)|=|RAI_{2n-1}(4321)|+|RAI_{2n-4}(4321)|.$

\endproof

The preceding lemma implies that the sequence $(|AI_{2n-1}(4321)|)_{n\geq 1}$ coincides with sequence \seqnum{A102071} in \cite{Sl}.

\begin{thm}\label{robs}
$$|AI_{n}(1234)|=|RAI_{n}(4321)|$$ and $$|RAI_{n}(1234)|=|AI_{n}(4321)|.$$
\end{thm}
\proof
The \textit{Robinson-Schensted} map $RS$ (see e.g. \cite{vL}) associates an involution $\pi$ with a standard Young Tableau whose descent set $D$ is equal to the descent set of $\pi.$ If we apply the inverse of the map  $RS$ to the transposed Tableau we get another involution $\pi'$ whose descent set is $[n-1]\setminus D.$ In particular, $\pi$ is alternating if and only if $\pi'$ is reverse alternating. Moreover, by the properties of the map $RS,$ $\pi$ avoids $1234$ if and only if the corresponding tableau has at most three columns, hence the transposed tableau has at most three rows and $\pi'$ avoids $4321.$
\endproof

Theorems \ref{RAI4321even}, \ref{RAI4321} and \ref{robs} and Lemma \ref{rc} imply the following result.

\begin{cor}\label{enum1234}
$$|RAI_{2n}(4321)|=|AI_{2n}(1234)|=M_n$$ and 
$$|RAI_{2n-1}(1234)|=|AI_{2n-1}(4321)|=|RAI_{2n-1}(4321)|=$$ $$|AI_{2n-1}(1234)|=M_n-M_{n-2}.$$
\end{cor}

We now give an explicit formula for the cardinalities of the sets $RAI_{2n}(1234)$ and $AI_{2n}(4321).$

\begin{thm}
$|RAI_{2n}(1234)|=|AI_{2n}(4321)|=M_{n+1}-2M_{n-1}+M_{n-3}.$
\end{thm}
\proof
We proceed as in the proof of Theorem \ref{RAI4321}.
Consider a permutation $\pi\in AI_{2n+1}(4321).$ The symbol $2n$ is either at position $2n$ or $2n+1,$ otherwise $\pi$ would contain $4321.$

If $\pi_{2n+1}=2n,$ the permutation $\pi'$ obtained from 
$\pi$ by removing the element $2n+1$ and standardizing is an arbitrary element of $AI_{2n}(4321)$ which ends with its maximum element.  

If $\pi_{2n}=2n,$ set 
$\pi=\left(\begin{smallmatrix} 
\cdots & k & \cdots & 2n-1 & 2n & 2n+1\\
\cdots & 2n+1 & \cdots & j & 2n & k
\end{smallmatrix}\right)$.

If $j<k,$ the permutation $\pi'$ obtained from 
$\pi$ by removing the element $2n$ and standardizing is an arbitrary element of $AI_{2n}(4321)$ which does not end with its maximum element. 

If $j>k,$ the element $j=2n-1$ must be fixed and $k=2n-2,$ otherwise $\pi$ would contain an occurrence of $4321.$ Removing from $\pi$ the last four elements  we get an arbitrary permutations in $AI_{2n-3}(4321).$

Hence $|AI_{2n+1}(4321)|=|AI_{2n}(4321)|+|AI_{2n-3}(4321)|$ and the assertion  now follows from Corollary \ref{enum1234}.
\endproof

\section{Pattern 3412}

\begin{thm}
$$|RAI_{2n}(3412)|=|AI_{2n+2}(3412)|=|AI_{2n+1}(3412)|=|RAI_{2n+1}(3412)|=M_n.$$
\end{thm}
\proof
We recall that also the set $I_n(3412)$ corresponds bijectively  to the set $\mathcal M_n$ via a map $\Psi$ whose definition coincides with the definition of the map $\Phi$ of Section \ref{sec4321} (see \cite{Ba6}). 
It is easily seen that  $\Psi$ maps the set $RAI_{2n}(3412)$ to the set of Motzkin paths in $\mathcal{M}_{2n}$ whose diods are either of the form $UH,$ $HD,$ or $UD.$ Hence the map $\widehat \Psi=\Delta \circ \Psi $ is a bijection between  $RAI_{2n}(3412)$ and $\mathcal{M}_{n},$
where $\Delta$ is the map defined in the proof of Theorem \ref{RAI4321even}.

We observe that any permutation in $AI_{t}(3412)$ begins with 1, hence, a permutation $\pi \in AI_{2n+1}(3412)$ can be written as $\pi=1\,\pi',$ where $\pi'$ is an arbitrary permutation in $RAI_{2n}(3412).$ Moreover, any permutation $\pi$ in $AI_{2n+2}(3412)$ ends with its maximum element, hence it can be can be written as $\pi=1\,\pi'\,2n+2$ where $\pi'$ is an arbitrary permutation in $RAI_{2n}(3412).$

The last equality follows from Lemma \ref{rc}.
\endproof

\section{Patterns 2143, 2134 and 1243}

\begin{thm}\label{enum1243one}
$$|AI_{2n}(1243)|=|AI_{2n}(2143)|=|AI_{2n}(2134)|=M_n$$ and $$|AI_{2n-1}(2134)|=|RAI_{2n-1}(1243)|=M_n-M_{n-2}.$$
\end{thm}
\proof
It is an immediate consequence of Lemma \ref{rc}, Lemma \ref{BWX} and Corollary \ref{enum1234}. 
\endproof

\begin{thm}
$$|AI_{2n+1}(1243)|=|AI_{2n+1}(2143)|=|RAI_{2n+1}(2143)|=|RAI_{2n+1}(2134)|=M_n.$$
\end{thm}
\proof

By Lemmas \ref{rc}, \ref{BWX} and \ref{Bonabij}, it suffices to show that $|AI_{2n+1}(1234)\cap \{\sigma \in S_{2n+1} \,|\, \sigma(2)=2n+1\}|=M_n.$

If $\pi$ is a permutation in this last set, removing $2n+1$ and $2$ from $\pi$ and standardizing we get an arbitrary permutation $\widehat{\pi}$ in $I_{2n-1}(1234)$ with either $\Des(\widehat{\pi})= \{3,5,7,\ldots 2n-3\}$ or $\Des(\widehat{\pi})= \{1,3,5,7,\ldots 2n-3\}.$ In other words, this permutation is reverse alternating or it is reverse alternating only to the right of the second position. Hence we have
$$|AI_{2n+1}(1234)\cap \{\sigma \in S_{2n+1} \,|\, \sigma(2)=2n+1\}|=$$ $$|I_{2n-1}(1234)\cap \{\sigma \in S_{2n-1} \,|\, \Des(\sigma)= \{3,5,7,\ldots 2n-3\}\}|+|RAI_{2n-1}(1234)|.$$

We want to show that the sum on the right hand side of the previous equality is equal to $M_n.$ We proceed by induction on $n$. The assertion is trivial when $n=2$ or $n=3.$ Assume that the assertion is true for every $m<n.$
We know by Corollary \ref{enum1234} that $|RAI_{2n-1}(1234)|=M_n-M_{n-2}.$
Consider a permutation $\rho$ in $I_{2n-1}(1234)\cap \{\sigma \in S_{2n-1} \,|\, \Des(\sigma)= \{3,5,7,\ldots 2n-3\}\}.$ Notice that  $\rho(2)=2n-2$ and $\rho(3)=2n-1,$ otherwise the permutation would contain the subsequence $\rho(1)\;\,\rho(2)\;\, 2n-2\;\,2n-1 \sim 1234.$ Removing $2,3,2n-2$ and $2n-1$ and standardizing we get an arbitrary permutation in 
$I_{2n-5}(1234)\cap \{\sigma \in S_{2n-5} \,|\, \Des(\sigma)= \{3,5,7,\ldots 2n-7\}\}$ or in $RAI_{2n-5}(1234).$ By the inductive hypothesis we have 
$$|I_{2n-5}(1234)\cup \{\sigma \in S_{2n-5} \,|\, \Des(\sigma)= \{3,5,7,\ldots 2n-7\}\}|+$$ $$|RAI_{2n-5}(1234)|=M_{n-2}.$$ This completes the proof. 

\endproof 

\begin{example}
Consider $\pi=5\,11\,9\,10\,1\,7\,6\,8\,3\,4\,2\in AI_{11}(1234).$
If we remove 11 and 2 from $\pi$ and standardize we get the involution $$\widehat{\pi}=4\,8\,9\,1\,6\,5\,7\,2\,3$$ with $\Des(\widehat{\pi})=\{3,5,7\}.$ Removing from $\widehat{\pi}$ the symbols 8,9,2 and 3 we get the involution $2\,1\,4\,3\,5\in RAI_5(1234).$

Consider now the permutation $\pi=9\,11\,7\,8\,5\,10\,3\,4\,1\,6\,2.$
If we remove 11 and 2 from $\pi$ and standardize we get the involution $$\widehat{\pi}=8\,6\,7\,4\,9\,2\,3\,1\,5\in RAI_9(1234).$$
\end{example}

\begin{thm}
$$|RAI_{2n}(2143)|=M_{n-1}.$$
\end{thm}
\proof
Let $\pi$ be a permutation in $RAI_{2n}(2143).$ The maximum of $\pi$ is in position $1,$ otherwise $2\,1\,\pi_{2n-1}\,\pi_{2n}$ would be an occurrence of $2143$ in $\pi.$ Write $\pi$ as $2n\,\widehat{\pi}\,1.$ Then the standardization of $\widehat{\pi}$ is an arbitrary permutation in $AI_{2n-2}(2143).$ 
The assertion now follows from Theorem \ref{enum1243one}.
\endproof

\section{Patterns 3421 and 4312}

Notice that the these two patterns are inverse of each other. Hence $AI_{n}(3421)=AI_n(4312)=AI_n(3421,4312)$ and the same is true for reverse alternating involutions. 
The following theorem shows that all these classes are enumerated by the Fibonacci numbers $F_n$ (sequence \seqnum{A000045} in \cite{Sl}).

\begin{thm}
$$|AI_n(3421,4312)|=|RAI_n(3421,4312)|=F_{n-1},$$ where $F_k$ is the $k$-th Fibonacci number. 
\end{thm}
\proof

Let $\pi \in RAI_{2n}(3421,4312).$ The last element of $\pi$ is odd and the position of $2n$ is odd by Lemma \ref{str}.

Notice that  the symbol $2n-1$ can be either at position $2n$ or $2n-1,$ otherwise $\pi$ would contain the subsequence $2n\;\,2n-1\;\,\pi(2n-2)\;\,\pi(2n-1)$ whose standardization is $4312.$

If $\pi_{2n}=2n-1,$ then $\pi=\pi'\;\,2n\;\,2n-1,$ where $\pi'$ is an arbitrary permutation in $RAI_{2n-2}(3421,4312).$

If $\pi_{2n-1}=2n-1,$ let $y=\pi_{2n-2}$ and $x=\pi_{2n}.$ If $x>y,$ removing  from $\pi$ the symbol $2n-1$ and standardizing we obtain an arbitrary permutation in $RAI_{2n-1}(3421,4312)$ which does not end with its maximum.

If $x<y,$ notice that either $\pi(2n-2)=2n-2$ or $\pi(2n-3)=2n-2$ (otherwise $\pi$ would contain the patterns). If $\pi(2n-2)=2n-2,$ then  $\pi(2n-3)$ is forced to be $2n$ (because $\pi$ is alternating),
so we can write $\pi=\tau\;\,2n\;\,2n-2\;\,2n-1\;\,2n-3,$ where $\tau\in RAI_{2n-4}(3421,4312).$ We  associate to $\pi$ the permutation $\pi'=\tau\;\, 2n-2\;\,2n-3\;\,2n-1.$ Such $\pi'$ is an arbitrary permutation in $RAI_{2n-1}(3421,4312)$ which ends with its maximum and in which $2n-3$ precedes the maximum.

If $\pi(2n-3)=2n-2,$ we  write
$\pi=\sigma' \;\, 2n \;\, \sigma'' \;\, 2n-2 \;\, 2n-3 \;\, 2n-1 \;\, x$
 and we  associate to $\pi$ the permutation $\pi'= \sigma' \;\, 2n-2 \;\, \sigma'' \;\, 2n-3 \;\, x \;\, 2n-1 $ which is an arbitrary permutation in $RAI_{2n-1}(3421,4312)$ which ends with its maximum and in which $2n-3$ follows the maximum.

As a consequence $$|RAI_{2n}(3421,4312)|=|RAI_{2n-1}(3421,4312)|+|RAI_{2n-2}(3421,4312)|.$$

Let $\pi \in RAI_{2n+1}(3421,4312).$ The last element of $\pi$ is odd and the position of $2n+1$ is odd by Lemma \ref{str}.
Notice that the symbol $2n+1$ is either at position $2n+1$ or $2n-1,$ otherwise $\pi$ would contain the subsequence $2n\;\,2n+1\;\,\pi(2n-1)\;\,\pi(2n)\sim 3421.$

If $\pi_{2n+1}=2n+1,$ then $\pi=\pi'\,2n+1,$ where $\pi'$ is an arbitrary permutation in $RAI_{2n}(3421,4312).$

If $\pi_{2n-1}=2n+1,$ write $\pi=\tau \;\, 2n+1 \;\, x \;\, 2n-1.$
Then $\pi'=\tau \;\, 2n-1 \;\, x $ is an arbitrary permutation in $RAI_{2n}(3421,4312)$ with $2n-1$ as a fixed point. We observed above that such permutations are in bijection with $RAI_{2n-1}(3421,4312).$

As a consequence 
\begin{equation}\label{eq1}
|RAI_{2n+1}(3421,4312)|=|RAI_{2n}(3421,4312)|+|RAI_{2n-1}(3421,4312)|.
\end{equation}

Hence we can conclude that $$|RAI_{m}(3421,4312)|=F_{m-1}.$$

By Lemma \ref{rc}, we have also $|AI_{2n+1}(3421,4312)|=|RAI_{2n+1}(3421,4312)|.$

Consider now $\pi\in AI_{2n}(3421,4312).$ By Lemma \ref{str}, $\pi_1$ is odd and the symbol 1 is in odd position. Notice that either $\pi(1)=1$ or $\pi(3)=1$  otherwise $\pi$ would contain the subsequence $\pi_1 \, \pi_2\, 2\, 1,$ whose standardization is 3421.

If $\pi_1=1,$ the standardization of $\pi_2\ldots\pi_{2n}$ is an arbitrary permutation in $RAI_{2n-1}(3421,4312).$ If $\pi(3)=1,$ write $\pi=3\,\pi_2\,1\,\tau.$ The standardization of the word $\pi_2\,3\,\tau$ is an arbitrary permutation $\pi'$ in $RAI_{2n-1}(3421,4312)$ which fixes the symbol $2.$

Observe that, if $\pi\in RAI_{2n-1}(3421,4312),$ then $\pi(1)=2$ or $\pi(2)=2,$ namely  $RAI_{2n-1}(3421,4312)=A\cup B$ where
$A=RAI_{2n-1}(3421,4312)\cap \{\pi \in S_{2n-1}\,|\, \pi(2)=1\}$ and $B=RAI_{2n-1}(3421,4312)\cap \{\pi \in S_{2n-1}\,|\, \pi(2)=2\}.$

The set $A$ corresponds bijectively to the set $RAI_{2n-3}(3421,4312)$ by removing the firs two elements, hence the set $B$ corresponds bijectively to the set $RAI_{2n-2}(3421,4312)$  by equation (\ref{eq1}).

Hence, $$|AI_{2n}(3421,4312)|=|RAI_{2n-1}(3421,4312)|+|RAI_{2n-2}(3421,4312)|$$ and $$|AI_{m}(3421,4312)|=F_{m-1}.$$
\endproof

\begin{example}
Consider the set $RAI_6(3421,4312),$ whose elements are 
$$\alpha=2\,1\,6\,4\,5\,3\quad\quad \beta=2\,1\,4\,3\,6\,5 \quad \quad \gamma=4\,2\,6\,1\,5\,3$$
$$\delta=4\,2\,3\,1\,6\,5 \quad \quad \epsilon=6\,2\,4\,3\,5\,1.$$
\end{example}
The permutations $\beta$ and $\delta$ correspond to $2\,1\,4\,3$ and $4\,2\,3\,1,$ namely, the elements of the set $RAI_4(3421,4312),$ whereas the permutations $\alpha, \gamma$ and $\epsilon$ correspond to $2\,1\,4\,3\,5,$ $4\,2\,5\,1\,3$ and $5\,2\,3\,4\,1,$ respectively, namely, the elements of $RAI_5(3421,4312).$

\section{Patterns 2431, 4132, 3241 and 4213.}

Notice that the two patterns $2431$ and $4132$ are inverse of each other.  Hence $AI_{n}(2431)=AI_n(4132)=AI_n(2431,4132)$ and the same is true for reverse alternating involutions and for the pair of patterns $(3241,4213).$

\begin{thm}\label{theorem0}
$$|RAI_{2n}(2431,4132)|=|RAI_{2n}(3241,4213)|=|RAI_{2n+1}(3241,4213)|=$$ $$|AI_{2n+1}(2431,4132)|= 2^{n-1}.$$
\end{thm}
\proof
Consider a permutation $\pi\in RAI_{2n}(2431,4132).$ We want to show that $\pi=\tau' \, 2n\,\tau'',$ where $\tau'$ is a permutation of the symbols
$1,2,\ldots,i-1$ and $\tau''$ is a permutation of the symbols $i,\ldots,2n-1.$

Suppose on the contrary that there exists an element $b$ in  $\tau'$ such that $b>i,$ so that $\pi=\ldots b\ldots 2n\ldots i.$ Let $a=\pi_{2n-1}.$ Since $\pi$ is alternating, we have $a>i,$ hence $2n-1$ follows $2n$ in $\pi$ and the subsequence $b\, 2n\, 2n-1\, i$ is order-isomorphic to 2431.

As a consequence, the word  $2n\,\tau''$ is itself a reverse alternating involution avoiding 132, and there exists only one of such permutations, as observed in Section \ref{ltr}, while the word $\tau'$ is an arbitrary permutation in $RAI_{j}(2431,4132),$ $j$ even.

We can conclude that $$|RAI_{2n}(2431,4132)|=\sum_{j \mbox{ even, }j<2n}|RAI_{j}(2431,4132)|=2^{n-1},$$
where the last equality follows from the previous one by induction. 

The facts that $$|RAI_{2n}(2431,4132)|=|RAI_{2n}(3241,4213)|$$ and $$|RAI_{2n+1}(3241,4213)|=|AI_{2n+1}(2431,4132)|$$ follow by Lemma \ref{rc}.

Consider now a permutation $\pi$ in $RAI_{2n+1}(4213,3241).$

Write $\pi=m_1\, w_1\,m_2\, w_2\ldots m_k\, w_k,$  where the $m_i's$ are the ltr maxima of $\pi$ and the $w_i's$ are words. Since $\pi$ is alternating and avoids the patterns 4213 and 3241, it follows that for every $i=1,2\ldots, k-1$ the word $w_i$ with  is non empty and that $w_i<w_{i+1},$ namely, every symbol in $w_i$ is smaller than any symbol in $w_{i+1}.$ 

As a consequence $m_iw_i$ is an interval for every $i.$
Obviously $m_k=2n+1$ and its position is $\pi_{2n+1}$ since $\pi$ is an involution. Suppose that $\pi_{2n+1}\neq 2n+1.$ Since $\pi_{2n}<\pi_{2n+1},$  the symbol $2n$ should appear to the left of $2n+1$ in $\pi,$ but this contradicts the fact that $m_kw_k$ is an interval. Hence $2n+1$ is the last element of $\pi$ and removing it from $\pi$ we get an arbitrary permutation in $RAI_{2n}(4213,3241),$ so we have
$$|RAI_{2n}(4213,3241)|=|RAI_{2n+1}(4213,3241)|.$$

\endproof

\begin{thm}\label{theorem1}
$$|RAI_{2n+1}(2431,4132)|=|AI_{2n+1}(3241,4213)| =\left\lfloor \dfrac{2^{n-1}\cdot 5}{3} \right\rfloor.$$
\end{thm}

To  prove  this theorem we  need the following lemma.
\begin{lem}\label{Lemma1}
The number of connected permutations in $RAI_{2l+1}(4132,2431)$ is $$\left \lfloor \frac{2l+1}{4}\right \rfloor.$$
\end{lem}
\proof
Let $\tau\in RAI_{2l+1}(4132,2431)$ be a connected permutations.
Consider the position $i$ of $2l+1$ in  $\tau.$ This position is odd by lemma \ref{str} and $i\neq 2l+1$ since $\tau$ is connected. 

If $i=2l-1,$ we can remove from $\tau$ the symbols $2l+1$ and $2l-1$ obtaining an arbitrary connected permutation $\tau$ in $RAI_{2l-1}(2431,4132).$

If $i<2l-1,$ let $a=\tau_{2l},$ $b=\tau_{2l-1}$ and $c=\tau_{2l-2}.$

We can write $\tau=\left(\begin{smallmatrix} 
\cdots & a & \cdots & i & \cdots & b & \cdots & 2l-2 & 2l-1 & 2l &2l+1\\
w_1 & 2l & w_2 & 2l+1 & w_3 & 2l-1 & w_4 &  c & b & a & i
\end{smallmatrix}\right)$ where, since $\tau$ is alternating, $c<b>a<i.$

If $bai\sim 312,$ since $\tau$ avoids 2431, the elements of $w_1$ and $w_2$ are smaller than the elements of $w_4\,c\,b\,a\,i$ that, in turn, are smaller than the elements of $w_3.$ But now, if $a\neq i-1,$ we would have an occurence of 4132. Hence $a=i-1,$ and this is impossible since $\tau$ is connected ($w_1$ would be a connected component of $\tau$).

Hence $bai\sim 213.$ Notice that we must have $cba\sim 231$ and we can write 
$\tau=\left(\begin{smallmatrix} 
\cdots & a & \cdots & c & \cdots & b  &\cdots & i &  \cdots & 2l-2 & 2l-1 & 2l &2l+1\\
w_1 & 2l & w_{2_1} & 2l-2 & w_{2_2} & 2l-1 & w_3 & 2l+1 & w_4 &  c & b & a & i
\end{smallmatrix}\right)$

where, as above, $w_1<w_4<w_3<w_{2_2}<w_{2_1}.$ 

The word $w_1$ is empty, since $\tau$ is connected. 
If $i\neq 2l-3,$ we have $\tau_{2l-3}>\tau_{2l-2},$
since $\tau$ is alternating. Hence $2l-3$ must follow $2l-2.$ So we can conclude that also $w_{2_1}$ is empty and $c=2.$ 
Since $\tau_2<\tau_3,$ 3 must follows 2 in $\tau$ so $b=3$ and $w_{2_2}$ is also empty.

Iterating these arguments one can show that
$$\tau= 2l\;2l-2\;2l-1\;2l-4\;2l-3\;\ldots l-1\;2l+1\; l-3\;l-2\;\ldots 2\;3\;1\;l+1.$$ In particular $i=l+1.$ Since $i$ is odd, this implies that $2l+1\equiv 1\,\mbox{ mod }4.$

As a consequence, the number of connected permutations in $RAI_{2l+1}(4132,2431)$ is the same as the number of  connected permutations in $RAI_{2l-1}(4132,2431)$ if $2l+1\equiv 3\mbox{ mod }4$ and increases by one if $2l+1\equiv 1\mbox{ mod }4.$

As an example, when $2l+1=21$ such permutations are

$$4\; 2\; 6\; 1\; 8\; 3\; 10\; 5\; 12\; 7\; 14\; 9\; 16\; 11\; 18\; 13\; 20\; 15\; 21\; 17\; 19$$ 
$$8\; 6\; 7\; 4\; 10\; 2\; 3\; 1\; 12\; 5\; 14\; 9\; 16\; 11\; 18\; 13\; 20\; 15\; 21\; 17\; 19$$
$$12\; 10\; 11\; 8\; 9\; 6\; 14\; 4\; 5\; 2\; 3\; 1\; 16\; 7\; 18\; 13\; 20\; 15\; 21\; 17\; 19$$
$$16\; 14\; 15\; 12\; 13\; 10\; 11\; 8\; 18\; 6\; 7\; 4\; 5\; 2\; 3\; 1\; 20\; 9\; 21\; 17\; 19$$
$$20\; 18\; 19\; 16\; 17\; 14\; 15\; 12\; 13\; 10\; 21\; 8\; 9\; 6\; 7\; 4\; 5\; 2\; 3\; 1\; 11$$

(notice that all of them but the last one have 21 at position 19 and the last one has 21 in position 11). When $2l+1=23$ such permutations are

$$4\; 2\; 6\; 1\; 8\; 3\; 10\; 5\; 12\; 7\; 14\; 9\; 16\; 11\; 18\; 13\; 20\; 15\; 22\; 17\; 23\; 19\; 21$$
$$8\; 6\; 7\; 4\; 10\; 2\; 3\; 1\; 12\; 5\; 14\; 9\; 16\; 11\; 18\; 13\; 20\; 15\; 22\; 17\; 23\; 19\; 21$$
$$12\; 10\; 11\; 8\; 9\; 6\; 14\; 4\; 5\; 2\; 3\; 1\; 16\; 7\; 18\; 13\; 20\; 15\; 22\; 17\; 23\; 19\; 21$$
$$16\; 14\; 15\; 12\; 13\; 10\; 11\; 8\; 18\; 6\; 7\; 4\; 5\; 2\; 3\; 1\; 20\; 9\; 22\; 17\; 23\; 19\; 21$$
$$20\; 18\; 19\; 16\; 17\; 14\; 15\; 12\; 13\; 10\; 22\; 8\; 9\; 6\; 7\; 4\; 5\; 2\; 3\; 1\; 23\; 11\; 21$$

(all of them have the symbol 23 in position 21).

Since when $2l+1=5$ there is only one of such permutations, the number of connected permutations in $RAI_{2l+1}(4132,2431)$ is $\left \lfloor \frac{2l+1}{4}\right \rfloor.$
\endproof

Now we turn to the proof of Theorem \ref{theorem1}.

\vspace{0.3 cm}
\noindent
\textit{Proof of Theorem }\ref{theorem1}.
Let $\pi\in RAI_{2n+1}(2431,4132).$ Consider the last connected component  $\tau$ of $\pi$ and write $\pi=\sigma\tau.$

Now, $\sigma$ is an arbitrary permutation in $RAI_{2k}(2431,4132),$ with $k\geq 0$ (notice that $\sigma$ has even length, since $\tau$ is connected and $\pi$ is alternating) and $\tau',$ the standardization of $\tau,$ is an arbitrary connected permutation in $RAI_{2n+1-2k}(2431,4132),$ with $l\geq 0.$

Hence, by Lemma \ref{Lemma1} and by Theorem \ref{theorem0}, we have 
$$|RAI_{2n+1}(2431,4132)|=2^{n-1}+\sum_{k\geq 1}^{n-1}2^{k-1}\left\lfloor \frac{2n+1-2k}{4}  \right \rfloor+\left\lfloor\frac{2n+1}{4} \right \rfloor, $$
where the first summand of the right hand side corresponds to the case $\tau=2n+1$ and the third summand corresponds to the case $|\tau|=2n+1.$

Now it can be shown by induction that $$|RAI_{2n+1}(2431,4132)|=\left \lfloor \frac{2^{n-1}\cdot 5}{3}\right \rfloor.$$

The fact that $|AI_{2n+1}(3241,4213)|=|RAI_{2n+1}(2431,4132)|$ follows by Lemma \ref{rc}.

\endproof

Notice that $\{|RAI_{2n+1}(2431,4132)|\}_{n}$ is sequence \seqnum{A081254} in \cite{Sl}.

On the contrary, the sequence enumerating $AI_{2n}(3241,4213)$ and  $AI_{2n}(2431,4132),$ whose first terms are $1,2,5,9,17,31,59,\ldots,$ is not present in \cite{Sl}.

\section{Patterns 2413 and 3142}

The set $AI_n(2413,3142)$ coincides with the set of alternating Baxter involutions of length $n.$ 

We recall that a \textit{Baxter permutation} of length $n$ is a permutation $\pi\in S_n$ such that, for every $1\leq i\leq j\leq k\leq  l\leq n,$ 
$$\mbox{if }\pi_i+1=\pi_l\mbox{ and }\pi_j>\pi_l\mbox{ then }\pi_k>\pi_l\mbox{ and}$$
$$\mbox{if }\pi_l+1=\pi_i\mbox{ and }\pi_k>\pi_i\mbox{ then }\pi_j>\pi_i.$$

In \cite{Ouc} the author shows that the set of doubly alternating Baxter permutations coincides with the set of doubly alternating permutations avoiding 2413 and $3142=2413^{-1}.$

An involution is clearly doubly alternating. As a consequence, the set $AI_n(2413,3142)$ coincides with the set of alternating Baxter involutions of length $n.$  A recurrence relation for this sequence has been found in \cite{Min} and the corresponding sequence in \cite{Sl} is  \seqnum{A347546}.

Moreover we have $$|AI_{2n+1}(2413,3142)|=|RAI_{2n+1}(2413,3142)|=|RAI_{2n}(2413,3142)|,$$ where the first equality follows from Lemma \ref{rc} and the second one is a consequence of the fact that each permutation in $RAI_{2n+1}(2413,3142)$ ends with its maximum (see \cite{Min}) and this maximum can be removed, hence obtaining any permutation in $RAI_{2n}(2413,3142).$

\section{Patterns 4123, 2341}

Since $4123=2341^{-1}$ we have $AI_n(4123)=AI_n(2341)=AI_n(4123, 2341)$ and the same is true for reverse alternating involutions.

To enumerate these classes we need the following lemma.

\begin{lem}
The number of connected permutations in $RAI_{n}(4123, 2341)$ is $\begin{cases} 1 & \mbox{ if }n\mbox{ is odd and }n\neq 3\\
2 & \mbox{ if }n\mbox{ is even and }n\geq 6\\
0 & \mbox{ if }n=3\\
1 & \mbox{ if }n=2,4.
\end{cases}.$\\
The number of connected permutations in $AI_n(4123,2341)$ is
$\begin{cases} 1 & \mbox{ if }n\geq 4\mbox{ or }n=1\\
0 & \mbox{ if }n=2\mbox{ or }3\\
\end{cases}.$
\end{lem}
\proof
The small cases are trivial, so we can assume $n\geq 5.$

Let $\pi\in RAI_{2n}(4123, 2341)$ be connected. Write $\pi=w_k\,m_k\, w_{k-1} m_{k-1}\ldots m_2\,w_1\,m_1$ where the $m_i's$ are the rtl maxima of $\pi.$

Since $\pi$ is reverse alternating and avoids 4123, every $w_i$ with $2\leq i\leq k-1$ must have length 1.
For the same reason, $w_1$ must be empty ($\pi$ ends with a descent) and the length of $w_k$ is 0 or 2.

If $|w_k|=0,$ $$\pi=2n\,2n-2\,2n-1\,2n-4\,2n-3\ldots 4\,5\,2\,3\,1, $$
if $|w_k|=2,$
$$\pi=2n-2\,2n-4\,2n\,2n-6\,2n-1\,2n-8\,2n-3\ldots 2\,7\,1\,5\,3.$$

In a similar way it is possible to prove that the only connected permutation in $RAI_{2n+1}(4123,2341)$ is 
$$\pi=2n\,2n-2\,2n+1\,2n-4\,2n-1\,2n-6\,2n-3\ldots 9\,4\,7\,2\,5\,1\,3. $$

The classification of the connected permutations in $AI_n(4123,2341)$ is fully analogous.  
\endproof

\begin{example}
There are two connected permutations in $RAI_{10}(4123,2341),$ namely, $$8\, 6\, 10\, 4\, 9\, 2\, 7\, 1\, 5\, 3 \mbox{ and } 10\, 8\, 9\, 6\, 7\, 4\, 5\, 2\, 3\, 1,$$ while the only connected permutation in $RAI_{11}(4123,2341)$ is $$10\, 8\, 11\, 6\, 9\, 4\, 7\, 2\, 5\, 1\, 3.$$
The only connected permutation in $AI_{10}(4123,2341)$ is 
$$9\, 10\, 7\, 8\, 5\, 6\, 3\, 4\, 1\, 2$$ and the only connected permutation in $AI_{11}(4123,2341)$ is $$9\, 11\, 7\, 10\, 5\, 8\, 3\, 6\, 1\, 4\, 2.$$
\end{example}

In the following theorem we find the ordinary generating functions that  enumerate reverse alternating involutions avoiding 4123 and 2341 of even and odd length (sequences \seqnum{A052980} and \seqnum{A193641} in \cite{Sl}, respectively). We also find the same generating function for the alternating even case. This last sequence does not appear in \cite{Sl}.

\begin{thm}
We have $$\sum_{n\geq 0 }|AI_{2n+1}(4123, 2341)|x^{2n+1}=\sum_{n\geq 0 }|RAI_{2n+1}(4123, 2341)|x^{2n+1}=$$ $$\frac{x^5-x^3+x}{1-2x^2-x^6}, $$
$$ \sum_{n\geq 0 }|RAI_{2n}(4123, 2341)|x^{2n}=\frac{1-x^2}{1-2x^2-x^6},$$
and $$\sum_{n\geq 0 }|AI_{2n}(4123, 2341)|x^{2n}=1+\frac{x^4}{1-x^2}+\frac{(x^5-x^3+x)^2}{(1-2x^2-x^6)\cdot (1-x^2)}.$$
\end{thm}

\proof
The fact that $|RAI_{2n+1}(4123, 2341)|=|AI_{2n+1}(4123, 2341)|$
follows from Lemma \ref{rc}.

Let $\pi \in RAI_{2n+1}(4123, 2341).$ Decompose $\pi$ in connected components as $\pi=\tau_1\tau_2\ldots \tau_k.$ 
Denote by $\widehat \tau_i$ the standardization of $\tau_i.$ The $\widehat \tau_i's$ are arbitrary connected permutations with $\widehat \tau_i\in \cup_{j\geq 1} RAI_{2j}(4123, 2341)$ for $1\leq i\leq k-1$ (because $\pi$ is reverse alternating) and $\tau_k\in \cup_{j\geq 0} RAI_{2j+1}(4123, 2341).$ By the previous lemma it follows that the ordinary generating function that counts connected permutations  in $\cup_{j\geq 1} RAI_{2j}(4123, 2341)$ (by length) is $F(x)=x^2+x^4+\frac{2x^6}{1-x^2}$ and the ordinary generating function that counts connected permutations  in $RAI_{2j+1}(4123, 2341)$ is $G(x)=x+\frac{x^5}{1-x^2}.$

As a consequence $$\sum_{n\geq 0 }|RAI_{2n+1}(4123, 2341)|x^{2n+1}=\frac{G(x)}{1-F(x)}.$$ Trivial algebraic manipulations lead to the desired generating function. 

Similarly, given $\pi \in RAI_{2n}(4123, 2341)$ we can decompose $\pi$ in connected components as $\pi=\tau_1\tau_2\ldots \tau_k$ where $\widehat \tau_i\in \cup_{j\geq 1} RAI_{2j}(4123, 2341)$ for $1\leq i\leq k.$

As a consequence $$\sum_{n\geq 0 }|RAI_{2n+1}(4123, 2341)|x^{2n+1}=\frac{1}{1-F(x)}.$$ 

The decomposition in connected components of a permutation  $\pi \in AI_{2n}(4123, 2341)$ is more subtle than the above ones because these components are not necessarily only of even or odd length. However such a permutation $\pi$ can only be
\begin{itemize}
    \item the empty one,
    \item connected itself,
    \item of the form $\sigma\tau,$ where $\sigma$ is an arbitrary permutation in $\cup_{j\geq 0} AI_{2j+1}(4123, 2341)$ and $\tau$ is a connected permutation in $\cup_{j\geq 0} RAI_{2j+1}(4123, 2341).$
\end{itemize}

By the above results and by the previous lemma we get
$$\sum_{n\geq 0 }|AI_{2n}(4123, 2341)|x^{2n}=1+\frac{x^4}{1-x^2}+\frac{x^5-x^3+x}{1-2x^2-x^6}\cdot \frac{x^5-x^3+x}{1-x^2}$$ where the three summand on the left hand side correspond to the three dotted cases.

\endproof

\section{Other patterns}

The following three conjectures are based on numerical evidences. The first one covers the last open case about the patterns 1243 and 2134 and the second covers the patterns 1432 and 3214.

\begin{conj}\label{conjecture1}
$$|RAI_{2n}(1243)|=|RAI_{2n}(2134)|=M_n.$$
\end{conj}

\begin{conj}\label{conjecture2}
$$|AI_{2n}(1432)|=|AI_{2n}(3214)|=|RAI_{2n}(1432)|=|RAI_{2n}(3214)|=M_n,$$
$$|AI_{2n+1}(1432)|=|RAI_{2n+1}(3214)|=M_n,$$
$$|AI_{2n-1}(3214)|=|RAI_{2n-1}(1432)|=M_{n}-M_{n-2}.$$
\end{conj}

\begin{conj}
Let $\tau$ be any permutation of $\{4,\ldots,m\}, m\geq 4.$ Then 
$$|AI_n(123\tau)|=|AI_n(321\tau)|.$$
\end{conj}

A classical subject in pattern avoidance is the study of Wilf-equivalent patterns. Let $\sigma$ and $\tau$ two patterns and $P_n$ a given subset of $S_n.$ Then $\sigma$ and $\tau$ are said to be \textit{Wilf-equivalent} on $P_n$ if $|P_n(\sigma)|=|P_n(\tau)|.$

A proof of Conjectures \ref{conjecture1} and \ref{conjecture2} would conclude the Wilf-classification  of alternating and reverse alternating involutions  avoiding a pattern of length 4. In fact, trivial numerical experiments show that the only Wilf-equivalences among all the other patterns not mentioned in the paper  are the trivial ones given by the reverse complement map and the inverse map. 

\section*{Acknowledgments}

We thank the anonymous referee for his very detailed revision and his valuable suggestions.

\addcontentsline{toc}{section}{Bibliography}
\bibliographystyle{plain}

\end{document}